\documentclass[10pt,oneside]{amsart}
\usepackage{amsfonts}
\usepackage{geometry}
\usepackage{graphicx}
\usepackage{textcmds}
\usepackage{palatino}
\usepackage{url}
\usepackage{amsmath,amsthm,amsfonts,amscd, amssymb, epsf} 
\usepackage{geometry}
 \theoremstyle{plain}    
 \newtheorem{thm}{Theorem}[section]
 
 \numberwithin{equation}{section} %% Comment out for sequentially-numbered
 \theoremstyle{plain}
 \theoremstyle{plain}    
 \theoremstyle{definition}
 \newtheorem{defn}[thm]{Definition}
 \newtheorem{as}[thm]{Assumption} 
 \theoremstyle{plain}    
 \newtheorem{prop}[thm]{Proposition} 
  
 \newtheorem{lem}[thm]{Lemma}
 \newtheorem{cor}[thm]{Corollary}
 \newtheorem*{cor*}{Corollary}
 \newtheorem*{conj*}{Conjecture}
 \newtheorem*{thm*}{Theorem}
\newcommand{\bl}{\begin{lem}}
\newcommand{\el}{\end{lem}}

\newcommand{\bml}{\begin{multline}}
\newcommand{\eml}{\end{multline}}

\newcommand{\beq}{\begin{equation}}
\newcommand{\eeq}{\end{equation}}
\newcommand{\bp}{\begin{prop}}
\newcommand{\ep}{\end{prop}}
\newcommand{\bd}{\begin{defn}}
\newcommand{\ed}{\end{defn}}
\newcommand{\pf}{\begin{proof}}
\newcommand{\epf}{\end{proof}}

\newcommand{\fd}{finite-dimensional  }

\newcommand{\field}[1]{\ensuremath{\mathbb{#1}}}
\newcommand{\CC}{\field{C}}
\newcommand{\df}{\equiv}
\newcommand{\NN}{\field{N}}

\newcommand{\HHHH}{\field{H}}

\DeclareMathOperator*{\lto}{\mathnormal{o}(1)}

\DeclareMathOperator{\lsum}{\mathnormal{\sum}}
\DeclareMathOperator*{\psum}{ \lsum^\prime}
\DeclareMathOperator{\gip}{\Gamma ^ \prime _\infty   }
\DeclareMathOperator{\gi}{\Gamma_\infty   }

\DeclareMathOperator{\HH}{\HHHH^3}

\DeclareMathOperator{\prV}{\mathnormal{\mathbf{P_{\infty}}}}
\newcommand{\PP}{\field{P}}

\newcommand{\RR}{\field{R}}

\newcommand{\ZZ}{\field{Z}}

\newcommand{\F}{\mathcal{F}}

\newcommand{\scz}{\mathcal{S}}
\newcommand{\pg}{\mathcal{P}}

\newcommand{\D}{\mathcal{D}}
\newcommand{\hil}{\mathcal{H}}

\newcommand{\mC}{\mathcal{C}}

\newcommand{\smat}{\mathfrak{S}}

\DeclareMathOperator{\rep}{Rep(\Gamma,\mathnormal{V})}
\DeclareMathOperator{\PSL}{PSL}

\DeclareMathOperator{\SL}{SL}

\DeclareMathOperator{\vol}{vol}

\DeclareMathOperator{\R}{Re}
\DeclareMathOperator{\I}{Im}
\DeclareMathOperator{\pc}{PSL(2,\CC)}

\DeclareMathOperator{\LOX}{lox}
\DeclareMathOperator{\cuspi}{\mathcal{CE}}

\DeclareMathOperator{\tr}{tr}

\DeclareMathOperator{\en}{\mathnormal{\mathcal{E}(T)}}
\DeclareMathOperator{\oen}{\mathnormal{\left|\mathcal{E}(T) \right|}}

\DeclareMathOperator{\ren}{\mathnormal{\mathcal{E}(R)}}

\DeclareMathOperator{\cinf}{\mathnormal{\PP}}

\DeclareMathOperator{\hs}{\mathnormal{\hil(\Gamma,\chi)}}
\DeclareMathOperator{\lp}{\mathnormal{\Delta}}

\newcommand{\0}{\bf{0}}

\newcommand{\ra}{\rightarrow}

\newcommand{\p}{\mathfrak{p}}
\begin{document}
%\sffamily
\title[Regularized determinants of the laplacian]{Regularized determinants of the laplacian for cofinite Kleinian groups with finite-dimensional unitary representations}
\author{Joshua S. Friedman}
\address{ Department of Mathematics and Sciences, United States Merchant Marine Academy, 300 Steamboat Road, Kings Point, NY  11024}
\email{CrownEagle@gmail.com}
\email{friedmanj@usmma.edu}
\email{joshua@math.sunysb.edu}
\maketitle
\thispagestyle{empty}
\begin{abstract}
For cofinite Kleinian groups (or equivalently, finite-volume three-dimensional hyperbolic orbifolds) with finite-dimensional unitary representations, we evaluate the regularized determinant of the Laplacian using  W. M\"{u}ller's regularization. We give an explicit formula relating the determinant to the Selberg zeta-function.   
\end{abstract}
%\tableofcontents{}
\section{Introduction}
The regularized determinant of the Laplacian has been well studied on Riemann surfaces. In the case of compact Riemann surfaces, D'Hoker and Phong \cite{Phong}, and Sarnak \cite{Sarnak} related the regularized determinant to the Selberg zeta-function.  

For non-cocompact cofinite Fuchsian groups (or equivalently, finite-area non-compact Riemann surfaces with elliptic fixed points)
Venkov, Kalinin, and Faddeev \cite{VenKalFad} defined a regularized determinant for the Laplacian $\lp$ and related the determinant to the Selberg zeta-function. They regularized the trace of the resolvent kernel using the theory of Krein's spectral shift function \cite{Krein1,BirmanKrein,Yafaev}.  

Efrat \cite{Efrat1,Efrat2} defined a regularized determinant for cofinite torsion-free Fuchsian groups with singular characters, and related it to the Selberg zeta-function.  His regularization was essentially based on the Selberg trace formula. Efrat's paper gave rise to an interesting question: Can the regularized determinant be defined \emph{cleanly} in terms of general operator theory? In the compact case, the answer is yes. Here zeta-regularization is defined in terms of the heat kernel, which is  of trace-class. In the non-compact case the heat kernel is not even Hilbert-Schmidt.

W. M\"{u}ller \cite{Muller,Muller3,Muller2,Muller1} applied Krein's theory to define a regularization of the determinant, a relative determinant $\det(H,H_{0})$ for two self-adjoint operators $H,H_{0},$ satisfying $\tr \left(e^{-Ht}-e^{-H_{0}}\right) < \infty.$ 

M\"{u}ller's regularization can be used for elliptic operators on non-compact manifolds. In  \cite{Muller1}, M\"{u}ller evaluates his determinant for the case of the Laplacian for  finite-area surfaces with hyperbolic ends (a class of surfaces that includes Riemann surfaces), and relates the determinant to Efrat's regularization (and hence to the Selberg zeta-function).

In \cite{Park} J. Park studies a closely related problem. He studies eta-invarients of Dirac operators, and relates the regularized determinant of the Dirac Laplacian to the Selberg zeta-function for odd-dimensional hyperbolic manifolds with cusps. Park also uses the regularized determinant to extract information about Selberg zeta-function. 

Regularized determinants have also been evaluated in the case of infinite volume Riemann surfaces, by Borthwick, Judge, and Perry \cite{Perry}.

In this paper, we evaluate M\"{u}ller's relative determinant of $\lp$  for the case of finite-volume three-dimensional hyperbolic orbifolds with finite-dimensional unitary representations. Or in other words the Laplacian acting on the Hilbert space of $\chi-$automorphic ($\chi$ is a finite-dimensional unitary representation) functions on hyperbolic three-space. We relate the determinant to the Selberg zeta-function using the appropriate version of the Selberg trace formula (proved previously in  \cite{Friedman1,Friedman2}).

We remark that zeta-regularization of determinants has  found application in quantum field theory, in the works of Dowker and Critchley \cite{DowkerCritchley}, Hawking \cite{Hawking}, Elizalde \emph{et al.} \cite{EORB}, and Bytsenko, Cognola and Zerbini \cite{Bytsenko}.

\subsection*{Main Results}
Next we define some of the basic notions needed to state our main results. A  Kleinian group  is a discrete subgroup of $\PSL(2,\CC) =   \SL(2,\CC) / \pm I. $ Each element of $\PSL(2,\CC) $ is identified with a M\"{o}bius transformation, and has a well-known action on hyperbolic three-space  $\HH$  and on its boundary at infinity \mdash the Riemann sphere $\PP^1$ (see \cite[Section 1.1]{Elstrodt})  . A Kleinian  group is \emph{cofinite} iff it has a fundamental domain $\F \subset \HH $ of finite hyperbolic volume. 

We use the following coordinate system for hyperbolic three-space, 
$\HH \df   \{(x,y,r)\in\RR^{3}~|~r>0 \} \df  \{ (z,r) ~| z \in \CC, ~r > 0 \} \df    \{z + rj\in\RR^{3}~|~r>0 \}, $ with the hyperbolic metric 
$$ ds^{2} \df \frac{dx^{2}+dy^{2}+dr^{2}}{r^{2}}, $$ and volume form 
$$ dv \df \frac{dx\, dy\, dz}{r^{3}}. $$
The Laplace-Beltrami operator is defined by
$$ \lp \df -r^{2}(\frac{\partial^{2}}{\partial x^{2}}+
\frac{\partial^{2}}{\partial y^{2}}+ \frac{\partial^{2}}{\partial r^{2}})+ r\frac{\partial}{\partial r}, $$ and it acts on the space of  smooth functions $f:\HH \mapsto V, $ where $V$ is a \fd complex vector space with inner-product $\langle~,~\rangle_V.$  

Suppose that $\Gamma $ is a cofinite Kleinian group and  $\chi \in \rep$ ($\rep$ is the space of \fd unitary representations of $\Gamma$ in $V$). Then the Hilbert space of \emph{$\chi-$automorphic} measurable functions is defined by  
\begin{multline*}
\hs  \df  \{ f: \HH \ra V ~|~ f(\gamma P) = \chi(\gamma) f(P)~\forall \gamma \in \Gamma, \\ P \in \HH, $ and $\left<f,f \right> \df \int_{\F} \left<f(P),f(P)\right>_V\,dv(P) < \infty \}. 
\end{multline*}
Here $\F$ is a fundamental domain for $\Gamma $ in $\HH$,  and $\left<~,~\right>_V $ is the inner product on $V.$ Finally, let   $ \lp = \lp(\Gamma,\chi) $ be the corresponding positive self-adjoint Laplace-Beltrami operator on $\hs.$

Next we briefly describe the motivation for the functional regularized determinant. 
Let $f(s)=\psum_{m \in \D} \lambda_m^{-s}$ be a sum over the non-zero eigenvalues of $\lp.$ Then \emph{formally} 
$$f'(0) = \left. \frac{d}{ds}f(s) \right|_{s=0} =\left.-\left( \sum_{m \in \D}' \log(\lambda) \lambda_m^{-s}  \right) \right|_{s=0}= -\left( \sum_{m \in \D}' \log(\lambda)  \right), $$
and $$e^{-f'(0)} = \prod_{\lambda_{m}\neq 0} \lambda_{m}. $$ 
With this \emph{formal} calculation in mind, one can think of $e^{-f^{\prime}(0)}$ as the regularized determinant. Now, typically, $f(0)$ does not even converge, but $f(s)$ does converge for $\R(s)$ sufficiently large. Analytic continuation gives a possible value for $f'(0).$

The formal argument above works well when the orbifold in question is compact. In the non-compact case we compare $\lp$ with another self adjoint operator, $\lp_{0},$  the self-adjoint extension of the operator 
$$\sum_{i=1}^{k_{\infty}} \left(-r^2 \frac{d^2}{dr^2}+r\frac{d}{dr}\right): \bigoplus_{i=1}^{k_{\infty}} C_0^\infty\left( [Y,\infty) \right) \mapsto \bigoplus_{i=1}^{k_{\infty}} L^2 \left([Y,\infty),r^{-3}dr \right) $$
with respect to Dirichlet boundary conditions ( $\{f \in C_0^\infty\left( [Y,\infty) \right)~|~f(Y)=0 \}$). See  \S\ref{ssLp0} for the definitions of the notations used above. 
%:curr0
Define the projection (onto the constant Fourier coefficient)  
$$\p_0:\hs \mapsto \bigoplus_{i=1}^{k_{\infty}} L^2 \left([Y,\infty),r^{-3}dr \right) $$ by
$$ \p_0[f](r) = \frac{1}{|\pg|}\int_{\pg} \prV f(x,y,r)~dx dy \quad \text{for}~ r \geq Y. $$
Once again, see  \S\ref{ssLp0} for the definitions of the notations used above.

The analogue of $f(s)$ for non-compact spaces, following M\"{u}ller, is the relative zeta-function 
$$
\zeta(s,\lp,\lp_{0}) \df \frac{1}{\Gamma(s)} \int_{0}^{\infty} t^{s-1}\left(\tr \left( e^{-\lp t}-e^{-\lp_0 t}\p_0   \right)-\dim \ker \lp \right)~dt.
$$
Here $\R(s)>2.$
Note that in order for the integral above to converge, we need to know the asymptotics of 
$$\tr \left( e^{-\lp t}-e^{-\lp_0 t}\p_0   \right)-\dim \ker \lp 
$$ 
at both $t=0$ and $t=\infty.$ These asymptotics are given in Lemma~\ref{lemGrowthTheta}.

Our main article of interest is the \emph{regularized characteristic polynomial,} $\det \left(\lp-(1-s^{2})\right),$ which we call the regularized determinant. For $\R(s)>2$ define  $$H(w,s)\df H(w,s,\lp,\lp_{0}) \df \frac{1}{\Gamma(w)} \int_{0}^{\infty} t^{w-1} \tr \left( e^{-\lp t}-e^{-\lp_0 t}\p_0\right) e^{t(1-s^{2})}~dt, $$ and following \cite{Sarnak}, we define  
$$
\det(\lp - (1-s^{2}))=e^{-\frac{\partial H}{\partial w}(0,s)}.
$$

Our main results are: 
\begin{thm*} Let $\Gamma$ be a cofinite Kleinian group with one cusp at infinity, and let $\chi \in \rep.$ Then there exists constants $C_{2},C_{3},D_{1},$ depending on $\Gamma$ and $\chi$ (they are explicitly determined in \S\ref{secSelZetaDet}) such that
\begin{multline*}
\log\det \left(\lp-(1-s^{2}) \right) = \log Z(s,\Gamma,\chi) + s \left( k(\Gamma,\chi)\log(Y) +C_{1}\right) \\ +\frac{l_\infty}{[\gi:\gip]}\log \Gamma(s+1) 
+ \Omega(s)  - \frac{C_{2}}{2}\log s-\frac{2}{3}C_{3}s^{3}-D_{1}.
\end{multline*}
Here $Z(s,\Gamma,\chi)$ is the Selberg zeta-function (see \S\ref{ssDefSZF}), $\Omega(s)$ is a meromorphic function (see Equation~\ref{eqOmEll}). The constant $Y$ comes about from the decomposition of $\F = \F_{Y} \cup \F^{Y}$ into a compact set $\F_{Y}$ and a noncompact cusp sector $\F^{Y}$ (see \S\ref{secNotation} for more details). 
\end{thm*}
The rest of the notation is defined in \S\ref{secSelZetaDet}.

\begin{cor*}
Let $\Gamma$ be a cofinite torsion-free Kleinian group with one cusp at infinity, and let $\chi \in \rep$ be a regular character. Then 
$$\det \left(\lp-(1-s^{2}) \right) = Z(s,\Gamma,\chi)\exp\left(-s^{3}\frac{\vol \left( \Gamma \setminus \HH \right)}{6\pi}+sL(\Lambda_{\infty},\psi) \right). $$
\end{cor*}
\noindent The constant $L(\Lambda_{\infty},\psi)$ comes about from regularity at a cusp, and its value is computed using KroneckerÕs second limit formula. It can be realized explicitly using the Siegel function  $g_{-v, u} \left( \tau \right), $ namely 
$$ L( \Lambda, \psi ) = \frac{-2 \pi}{y} \log \left| g_{-v, u} \left( \tau \right) \right|. $$ See \S\ref{sssSiegel} for more details.

\begin{cor*}
Let $\Gamma$ be a cocompact Kleinian group, $\chi \in \rep.$ Then 
$$\det \left(\lp-(1-s^{2}) \right) = Z(s,\Gamma,\chi)\exp\left(-s^{3}\frac{\vol \left( \Gamma \setminus \HH \right)}{6\pi}+sC_{E} \right), $$
where $$C_{E}=\sum_{ \{R \} \text{\emph{nce}}}\frac{\tr _{V}\chi (R) \log N(T_{0})}{4|\ren |\sin ^{2}(\frac{\pi k}{m(R)})}.$$
\end{cor*}
\noindent The constant $C_{E}$ is related to the non-cuspidal elliptic elements of $\Gamma.$ See \S\ref{sssNCE} for more details.

I would like to thank Professor Leon Takhtajan for originally suggesting this problem to me, for reading over this paper, and for useful discussions. I would also like to thank Professor Werner M\"{u}ller for answering some questions of mine related to the writing of this paper. 

I would also like to thank the anonymous referee, for pointing out some important errors, and for suggesting very useful comments.

\section{General Definition of the Relative Zeta-function}
In this section (following \cite{Muller}) we state some basic facts concerning Krein's spectral shift function, and show how they lead to the definition of the general relative zeta-function. Later, in \S\ref{secRelZeta}, we specialize to $\lp$ on $\hs.$  For more details on the spectral shift function see \cite{Krein1,BirmanKrein,Yafaev}.

We first establish some notation. For $B,$  a self-adjoint operator on $\mathcal{H},$  $\sigma(B)$ and $\sigma_{\text{ess}}(B)$ are the spectrum and essential spectrum of $B$ respectively. 

Let $A$ and $A_0$, be bounded self-adjoint operators in $\mathcal{H}.$ Suppose that $V \df A-A_0$ is of trace-class. Let $R_0(z)=\left( A_0-z \right)^{-1}$ be the resolvent $A_0.$ The \emph{spectral shift function} of $A$ and $A_0$ \beq \label{Essf}
\xi(\lambda)=\xi(\lambda;A,A_0)=\pi^{-1} \lim_{\epsilon \ra 0} \arg \det \left( 1+V R_0(\lambda+i \epsilon) \right),
\eeq
exists for a.e. $\lambda \in \RR,$ is real-valued, and belongs to $L^1(\RR).$ The determinant in \eqref{Essf} is the \emph{Fredholm determinant.} In addition, 
$$
\tr(A-A_0) = \int_\RR \xi(\lambda)~d\lambda, \quad \| \xi \| \leq \| A-A_0 \|_1, 
$$
where $\|\cdot \|_{1}$ is the trace norm. 

The theory of the spectral shift function is reminiscent of the Selberg trace formula.
Let 
$$  
\mathcal{G}=\left\{ \phi:\RR \ra \RR~|~\phi\in L^1~\text{and} \int_\RR |\widehat{\phi}(p)|(1+|p|)~dp < \infty \right\}. $$
Then for every $\phi \in \mathcal{G},~\phi(A)-\phi(A_0)$ is a trace class operator and 
$$ \tr( \phi(A)-\phi(A_0)) = \int_R \phi'(\lambda) \xi(\lambda) ~d\lambda. $$

Part~\ref{lemLoc93923} of the next lemma will allow us to define the relative zeta function for $H,H_{0}.$
\bl \cite[Page 315]{Muller}  \label{LtraceHeatGrowth} Let $H,H_0$ be two non-negative self-adjoint operators in $\mathcal{H}$ and assume that $e^{-tH} - e^{-tH_0}$ is a trace class operator for $t>0.$ Then there exists a unique real valued locally integrable function $\xi(\lambda)=\xi(\lambda;H,H_0)$ on $\RR$ such that for each $t>0,$ $e^{-t\lambda}\xi(\lambda) \in L^1(\RR)$ and the following conditions hold:
\begin{enumerate}
\item $\tr( e^{-tH}-e^{-tH_0} ) = -t \int_0^\infty e^{-t\lambda} \xi(\lambda)~d\lambda.$ 
\item For every $\phi \in \mathcal{G},~\phi(H)-\phi(H_0)$ is a trace class operator and 
$$\tr(\phi(H)-\phi(H_0) ) = \int_\RR \phi'(\lambda)\xi(\lambda)~d\lambda. $$
\item \label{lemLoc93923} In addition, suppose $\sigma_{ess}(H_0)\subset [c,\infty),$ where $c>0;$ then $\ker H$ and $\ker H_0$ are both finite-dimensional, and there exists $c_1 > 0$ such that 
$$\tr(e^{-tH}-e^{-tH_0}) = \dim \ker H - \dim \ker H_0 + O(e^{-c_1 t}) $$
as $t \ra \infty.$
\end{enumerate}
\el

Let $h= \dim \ker H - \dim \ker H_0. $ Then it follows from Lemma \ref{LtraceHeatGrowth} that for $\R(s)>0,$ the integral 
$$\int_0^\infty t^{s-1} \left( \tr \left(e^{-tH}-e^{-tH_0} \right) - h \right)~dt $$
converges absolutely.

\bd \cite[Page 317]{Muller}
Suppose that $\sigma_{ess}(H_0)\subset [c,\infty),$ where $c>0.$ Then for $\R(s) > 0$, the relative zeta-function of $H$ and $H_0,$ $\zeta(s;H,H_0)$ is defined by 
$$\zeta(s;H,H_0) = \frac{1}{\Gamma(s)}\int_0^\infty t^{s-1} \left( \tr \left(e^{-tH}-e^{-tH_0} \right) - h \right)~dt$$ 
\ed

\section{The Operators $\lp,~\lp_{0}$ and their Heat Kernels}
The main goal of this section is to give an explicit formula for 
$$\tr (e^{-\lp t}-e^{-\lp_{0} t}\p_{0}).$$ 

\subsection{Notation} \label{secNotation}
Before we can define $\lp_{0}$ we must establish some notation (see \cite{Friedman1,Friedman2} for more details). In order to simplify our notation (and to make our paper more readable) we present our results under the assumption:
\begin{as} \label{asOne}
The cofinite Kleinian group $\Gamma$ has only one class of cusps at $\zeta = \infty \in \PP,$ and $\chi \in \rep.$ The lattice associated with $\zeta=\infty$ is  
$$\Lambda_\infty = \ZZ \oplus \ZZ \tau_\alpha,~\I(\tau_\alpha) > 0.$$ 
\end{as}

Let  $\Gamma_\infty < \Gamma $ denote the stabilizer subgroup  of the cusp at infinity $\zeta=\infty,$ 
$$
\gi \df \{~ \gamma \in \Gamma ~ | ~ \gamma(\infty) = \infty ~\},
$$ 
and let $ \gip $ be the maximal torsion-free parabolic subgroup of $\gi. $ By definition (of a cusp), $\gip$ is a free abelian group of rank two. The possible values for the index of $[\gi:\gip]$ are 1,2,3,4, and 6.  See \cite{Elstrodt}.

The subgroup $\gip$ is canonically isomorphic to a lattice  $\Lambda_\infty = \ZZ \oplus \ZZ \tau_\alpha. $ Without loss of generality we can assume that $\I(\tau_\alpha) > 0. $ Let $\epsilon$ be root of unity of order $[\gi:\gip]$. Then 
\begin{enumerate}
\item $$  \gip =  
\left\{ \,\left.\left(\begin{array}{cc}
1 &  b\\
0 &   1\end{array}\right)\,\right|\, b \in \Lambda_{\infty} \,\right\},  $$
\item $$  \gi = \left\{ \,\left.\left(\begin{array}{cc}
\epsilon^{n} & \epsilon^{n} b\\
0 & \epsilon^{-n}\end{array}\right)\,
\right|\, b \in \Lambda_{\infty},n = 0,\dots [\gi:\gip] \right\} /\{\pm I\}.$$
\end{enumerate}

Let $\pg \subset \CC $ be a fundamental domain\footnote{The set $\pg$ is a euclidean polygon. \\} for the action of $\gi$ on $\CC,$ and let $\pg^\prime $ be the fundamental parallelogram with base point at the origin for the lattice $\Lambda_\infty. $ For $Y>0$ set 
$$
F^{Y} \df \F(Y) \df \{\, z+rj\,|\, z\in\pg,\, r\geq Y\,\}.$$ Then for $Y$ sufficiently large, there exists a compact set $\F_{Y},$ disjoint from $\F^{Y},$  so that $\F \df \F_{Y} \cup \F^{Y}$ is a fundamental domain for $\Gamma.$  

Define the \emph{singular} space by 
\beq \label{eqSingular}V_{\infty} \df \{ v \in V \, | \, \chi(\gamma)v=v, \, \, \, \forall \gamma \in \gi \,\}, \eeq  and the \emph{almost singular} space 
\beq \label{eqAlmostSingular} V_{\infty}^\prime \df \{ v \in V \, | \, \chi(\gamma)v=v, \, \, \, \forall \gamma \in \gip \,\}.  \eeq  

If $\dim V_{\infty} > 0,$ then $\chi$ is called \emph{singular} with index of singularity $k(\Gamma,\chi) \df k_{\infty} \df \dim V_{\infty}.$ If $\dim V_{\infty} = 0,$ then $\chi$ is called \emph{regular}.
Set $l_{\infty} \df \dim V'_{\infty}.$
Let $\prV$ denote the orthogonal projection $$\prV:V \mapsto V_{\infty}.  $$

Fix an orthonormal basis $ \{ v_i \}_{i=1}^{k_\infty} $ for $V_\infty. $ For $P \in \HH,\, \R(s)>1,$ and $i=1\dots k_{\infty};$ we define the  \emph{Eisenstein series}  by 
$$
E_{i}(P,s)   \df  E(P,s,i,\Gamma,\chi) \df  \sum_{M \in \gi \setminus \Gamma } \left(r(MP)
\right)^{1+s}\chi(M)^{*}v_i. $$

The series $E_{i}(P,s) $ converges uniformly and absolutely on compact subsets of  $\{ \R(s)> 1 \} \times \HH $ to a $\chi-$automorphic function  that  satisfies
$$ \lp E(~\cdot~,s,\alpha ,v) = \lambda E(~\cdot~,s,\alpha ,v), \quad \lambda=1-s^{2}, $$  and admits a meromorphic continuation to the whole complex plane \cite{Friedman1}.

For $P = z+rj,~P'=z'+r'j \in \HH$ set 
$$ \delta(P,P') \df  \frac{|z-z'|^{2}+r^{2}+ r'^{2}}{2rr'}. $$ It follows that $\delta(P,P') = \cosh(d(P,P'))$, where $d$ denotes the hyperbolic distance  in $\HH.$  Next, for $k \in \scz \df \scz([1,\infty))$ a Schwartz-class\footnote{The space of smooth functions $k:[1,\infty) \rightarrow \CC$ that satisfy $ \lim_{x \ra \infty} x^n k^{(m)}(x) = 0 $ for all $n,m \in \NN_{\geq 0}.$} function, define $K(P,Q)$ by $$K(P,Q) = k(\delta(P,Q)).$$ The function $K(P,Q)$ is called  a \emph{point-pair invariant.} 
Set
$$K_\Gamma(P,Q) \df \sum_{\gamma\in\Gamma} \chi(\gamma)K(P,\gamma Q).$$  
The decay properties of the function $k$ guarantee that the series above converges absolutely and uniformly on compact subsets of
$\HH \times \HH$ \cite[Theorem 6.4.1]{Elstrodt}.  The function $K_\Gamma(P,Q)$ is the kernel of a bounded operator $\mathcal{K} : \hs \mapsto \hs.$

The function $k$ leads to two other useful function: $h,$ the Selberg\ndash Harish-Chandra transform of $k;$ and $g,$ the Fourier transform of $h.$ Explicitly:  
\beq \label{eqSHC}
h(\lambda) = h(1-s^2) \df \frac{\pi}{s}
\int_{1}^{\infty}k\left(\frac{1}{2}\left(t+\frac{1}{t}\right)\right)
(t^{s}-t^{-s})\left(t-\frac{1}{t}\right)\,\frac{dt}{t},~~\lambda = 1-s^2,
\eeq
and for $r\in \RR$ set 
$$ g(r) = \frac{1}{2\pi} \int_{\RR} h(1+x^2)e^{-ixr}\,dx. $$  

For $v,w \in V $ let $ v \otimes \overline{w}$ be the linear operator in $V$ defined by    
$ v \otimes \overline{w}(x) = <x,w>v.  $
An immediate application of the Spectral Decomposition Theorem (\cite{Friedman1,Friedman2}), and the Selberg\ndash Harish-Chandra Transform yields (see \cite[Equation 6.4.10, page 278]{Elstrodt}):
\begin{lem} \label{lemKerExpan}
Let  $k \in \scz $ and  $h:\CC\ra\CC$ be the Selberg--Harish-Chandra  Transform
of $k.$ Then 
\begin{multline}
K_\Gamma(P,Q)   = \sum_{m \in \D} h(\lambda_m)e_m(P) \otimes 
\overline{e_m(Q)}   \\ +
\frac{1}{4\pi} \sum_{l = 1}^{k_\infty}
\frac{1}
{\left|\pg \right|}
\int_{\RR} h \left( 1+x^2 \right) E_{l}(P,ix) \otimes \overline{E_{l}(Q,ix)} \, dt, 
\end{multline} where $|\pg|$ denotes the euclidian area of $\pg \subset \CC.$
The sum and integrals converge absolutely and uniformly on compact subsets 
of $\HH \times \HH$.
\end{lem}

We conclude this section with some notation that will be needed  to state our main result.
\subsubsection{Regular representations and Siegel's theta function} \label{sssSiegel}
Recall that
$\Lambda_{\infty} =  \ZZ \oplus \ZZ \tau \subset \CC$ with $\I(\tau) > 0.$ It follows that $\chi$ restricted to (the abelian group) $\Lambda_{\infty}$ diagonalizes into  characters\footnote{One-dimensional unitary representations.} $\psi_{l}$ for $l=1\dots l_{\infty},$ and the identity character for $l=l_{\infty}+1 \dots \dim V.$ For each $\psi_{l},$  $u_{l},v_{l} \in \RR $ are not both integers, satisfying  
$  \psi_{l}(1) =  e^{2 \pi i u_{l}} ~~\text{and}  $
$ \psi_{l}(\tau) = e^{2 \pi i v_{l}}. $ We define 
$$ L( \Lambda_{\infty}, \psi_{l} ) = \frac{-2 \pi}{y} \log \left| g_{-v_{l}, u_{l}} \left( \tau \right) \right|, $$
where $g_{a_1,a_2} $ is the Siegel function, $$
g_{a_1,a_2}(\tau) = -q_{\tau}^{(1/2)\textbf{B}_{2} (a_1)} e^{ 2 \pi i a_2 (a_1-1)/2 }(1-q_z) \prod_{n=1}^{\infty} (1 - q_{\tau}^n q_z)(1 - q_{\tau}^n /q_z),
$$
$ \textbf{B}_{2} (X) = X^2 -X + 1/6, $ $ q_{\tau} = e^{ 2 \pi i \tau},$  $ q_z =  e^{ 2 \pi i z}, $ and $z = a_1\tau+a_2. $

\subsubsection{Non-cuspidal elliptic elements} \label{sssNCE}
Let $\{R\}_{\text{nce}}$ be a set of representatives of the non-cuspidal elliptic elements of $\Gamma,$ the elliptic elements that do not fix a cusp ($\infty$ under Assumption~\ref{asOne}). Following \cite[Definition 5.3.2]{Elstrodt}, the \emph{Elliptic number} of $\Gamma$ is
$$\sum_{ \{R \}\text{nce}}\frac{\tr _{V}\chi (R) \log N(T_{0})}{4|\ren |\sin ^{2}(\frac{\pi k}{m(R)})}.$$
For a fixed representative $R,$ $N(T_0)$ is the minimal norm of a hyperbolic or loxodromic element of the centralizer $\mathcal{C}(R).$  The element $R$ is understood to be a $k-$th power of a primitive non cuspidal elliptic element $R_0 \in \mathcal{C}(R)$ describing a hyperbolic rotation around the fixed axis of $R$ with minimal rotation angle $\frac{2 \pi}{m(R)}.$  Further, $\ren$ is the maximal finite subgroup contained in $  \mathcal{C}(R).$

\subsubsection{Cuspidal elliptic elements} \label{sssCE}
Denote by $ \cuspi $  set of elements of $\Gamma$ 
which are $\Gamma$-conjugate to an element of $ \Gamma_\infty
\setminus \Gamma_\infty^\prime = \{ \gamma \in \gi~|~\gamma~\text{is not parabolic nor the identity element}~\}. $  
We fix representatives  of conjugacy classes of  $\cuspi,~
g_{1}, \dots , g_{d}\footnote{There are only finitely many distinct conjugacy classes of elliptic elements in a cofinite Kleinian group.} $ that have the form 
\beq g_{i} =  \left(\begin{array}{cc}
 \epsilon_{i} & \epsilon_{i} \omega_{i} \\
 0 & \left( \epsilon_{i}  \right)^{-1}\end{array}
\right).
\eeq
Let $\mC (g) $ denote the centralizer in $\Gamma $ of an element $g \in \cuspi.$  In addition, let   $\{ p_i,\infty \} $ be the set of fixed points in $\cinf$ of the element $g_i.$  Since $g_i $ is a cuspidal elliptic element it follows that  $p_i $ is a cusp of $\Gamma $ (see \cite{Elstrodt} page 52).  Hence by Assumption~\ref{asOne} there is an element $\gamma_i \in \Gamma $ with $\gamma_i \infty = p_i. $  Let $c_i $ is the lower left hand (matrix) entry of $\gamma_i.$

\subsection{The Heat Kernel of $\lp$ as a Poincar\'{e} Series} The heat kernel for\footnote{We abuse notation and allow $\lp$ to represent both the self-adjoint operator on $\hs$ and the standard differential operator on smooth functions of $\HH.$}  $\lp$ on $\HH$ is a function $$u:\HH\times\HH\times (0,\infty) \mapsto \RR$$ satisfying 
$$e^{-\lp t}f(P) = \int_{\HH} u(P,Q,t) f(Q) ~dv(P), $$ for all $f$ in the domain of the self-adjoint operator $e^{-\lp t}.$ It is a classical result \cite{Davies1} that 
\beq \label{eqHeatKer}u(P,Q,t) = u(\rho,t)= (4\pi t)^{-3/2} \frac{\rho}{\sinh \rho} \exp \left( -t - \frac{\rho^2}{4t}    \right), \quad \text{where}~\rho = d(P,Q). \eeq 
In order to apply the theory of point pair invariants, we need to find a Schwartz-class function $k_t$ so that 
$$k_t(\delta(P,Q))=k_t(\cosh(d(P,Q))=u(P,Q,t).  $$ Using the formula:
$$\cosh^{-1}(x)=\ln(x+\sqrt{x^2-1}),\quad \text{for}~x \geq 1, $$ 
we obtain 
\beq \label{eqHeatPointPair} k_t(x)=\frac{e^{-t}}{(4\pi t)^{3/2}}\frac{\ln(x+\sqrt{x^2-1})}{\sqrt{x^2-1}}\exp\left(\frac{-\left(\ln(x+\sqrt{x^2-1}) \right)^2}{4t} \right).
\eeq
Observe\footnote{The limit follows from applying l'H\^{o}pital's rule to either \eqref{eqHeatKer} or \eqref{eqHeatPointPair}. The fact that the singularity cancels out is one advantage to working the the heat kernel instead of the resolvent kernel, which has a singularity and must be iterated. } that as $x \ra 1^{+}, k_t(x) \ra \frac{e^{-t}}{(4\pi t)^{3/2}},$ and that $k_t \in \scz.$ We have:
\bl
Let $$K_\Gamma(P,Q,t,\chi)=\sum_{\gamma \in \Gamma} \chi(\gamma) k_t(\delta(P,\gamma Q)).$$ Then $K_\Gamma(P,Q,t,\chi)$ is the heat kernel for $\lp$ on the Hilbert space $\hs.$
\el 
%:Spectral Decomposition of Heat Kernel

It follows from Equation~\ref{eqSHC} and Lemma~\ref{lemKerExpan} that 
$h(x)=e^{-t x}. $ By definition 
$$
g(r)=\frac{1}{2\pi} \int_\RR e^{-(1+x^2)t}e^{-ixr}~dx = \frac{\exp(-t)}{\sqrt{4\pi t}}\exp\left(-\frac{r^2}{4t}\right).
$$
By applying Lemma~\ref{lemKerExpan} we obtain the spectral expansion of  $K_{\Gamma}:$
\begin{multline}
\label{EheatKer}
K_\Gamma(P,Q,t,\chi)  = \sum_{m \in \D} e^{-\lambda_m t} e_m(P) \otimes 
\overline{e_m(Q)}   \\ +
\frac{1}{4\pi} \sum_{l = 1}^{k_{\infty}}
\frac{1}
{\left| \pg \right|}
\int_{\RR}  \exp\left( -(1+x^2)t \right) E_{l}(P,i x) \otimes \overline{E_{l}(Q,i x)} \, dx.
\end{multline} 

\subsection{The Operator $\lp_{0}$ and its Heat Kernel} \label{ssLp0}
For $\F \df \F_{Y} \cup \F^{Y},$ let $\lp_0$ be the self-adjoint extension of the operator 
$$\sum_{i=1}^{k_{\infty}} \left(-r^2 \frac{d^2}{dr^2}+r\frac{d}{dr}\right): \bigoplus_{i=1}^{k_{\infty}} C_0^\infty\left( [Y,\infty) \right) \mapsto \bigoplus_{i=1}^{k_{\infty}} L^2 \left([Y,\infty),r^{-3}dr \right) $$
with respect to Dirichlet boundary conditions ( $\{f \in C_0^\infty\left( [Y,\infty) \right)~|~f(Y)=0 \}$). Note that $\lp_{0}$ depends on $Y.$ It is understood that $\lp_{0}$ acts componentwise with respect to the basis for $V_{\infty}$ fixed in \S\ref{secNotation}.
The operator $\lp_0$ can be thought of as a Laplacian operator in its own right. In fact $\lp_0$ is closely related to the restriction of $\lp$ to $\F^{Y}$.  Define $$\p_0:\hs \mapsto \bigoplus_{i=1}^{k_{\infty}} L^2 \left([Y,\infty),r^{-3}dr \right) $$ by
$$ \p_0[f](r) = \frac{1}{|\pg|}\int_{\pg} \prV f(x,y,r)~dx dy \quad \text{for}~ r \geq Y. $$
\bl \label{lemInter}
For all $f$ in the domain of $\lp,$ $$\p_0 [\lp f] = \lp_0 \p_0[f].$$
\el
\pf
The proof follows from \cite[Pages 236-237]{Elstrodt} and the definition of $\p_0.$ 
\epf

For $\F = \F_Y\cup \F^{Y},~t > 0,~P=z+rj,P'=z'+r'j\in \F$ define  

\beq \label{eqSmallHeatKer} k(P,P',t) \df \left\{
\begin{array}{cc}
\prV  \frac{r r'\exp\left( -t \right)}{|\pg| \sqrt{4\pi t}}\left[ \exp\left(-\frac{\log^2(r/r')}{4t} \right) -\exp\left( -\frac{\left( \log(r r')-2\log(Y) \right)^2}{4t}   \right)\right]   & \mbox{for $P,P' \in \F^{Y},$}    \\ 
\\ 
 0 & \mbox{else}      
\end{array}
    \right. \eeq
The classical theory of the Heat Equation, for the half line, tells us that $k$ is the heat kernel of $\lp_{0}.$ In other words: 
\bl
For all $f \in \hs,$ 
$$e^{-\lp_0 t}\p_0[f] = \int_\F k(\,\cdot \,,P',t) f(P')~dv(P').  $$
\el

\subsection{The Regularized Trace} In this section we will prove that the \emph{regularized heat kernel}  $$ e^{-\lp t}-e^{-\lp_0 t}\p_0  $$ is of trace-class. Then using ideas from the standard proof of the Selberg trace formula, we will evaluate the trace explicitly. But before we can proceed, we need Theorem~\ref{thmHeatKernBounds}. Its proof is based on the classical Poisson Summation Formula, and it is   proved in Appendix~\ref{appHeatKernBounds}. 
\begin{thm} \label{thmHeatKernBounds}
For $P=z+rj, P'=z'+r'j \in \F^{Y},$ 
$$K_\Gamma(P,P',t,\chi)=\prV rr'\frac{\exp(-t)}{|\pg| \sqrt{4\pi t}}\exp\left(-\frac{\log^{2}(r/r')}{4t}\right)+O(1). $$ 
\end{thm}

\begin{thm} \label{thmTraceClass}
The operator $$e^{-\lp t}-e^{-\lp_0 t}\p_0$$ is of trace-class. \end{thm}
\begin{proof}
The proof is based on the decay properties of 
Theorem~\ref{thmHeatKernBounds} and Equation~\ref{eqSmallHeatKer}, and a clever trick, using the semi-group properties of the heat kernels, of Deift-Simon. We follow \cite[Page 259]{Muller3} and \cite[Prop. 2.1]{Park}. First note that by Theorem~\ref{thmHeatKernBounds} and Equation~\ref{eqSmallHeatKer}, $e^{-\lp t}-e^{-\lp_0 t}\p_0$ is Hilbert-Schmidt.

Next write $e^{-\lp t}-e^{-\lp_0 t}\p_0$ as 
\beq \label{eqTmp8nhgs}
e^{-\lp \tau}\left( e^{-\lp \tau}-e^{-\lp_0 \tau}\p_0 \right) + \left( e^{-\lp \tau}-e^{-\lp_0 \tau}\p_0 \right)e^{-\lp_{0} \tau}, 
\eeq
 where $\tau = t/2.$
Next choose a function $f \in C^{\infty}(\F)$ so that $0 < f \leq 1,$ $f(P)=1$ if $P \in \F_{Y},$ $F(z+rj)=(Y/r)^{1/4}$  if $P=z+rj \in \F^{Y}.$ Let $m_{f}$ be the multiplication operator by $f.$ Now rewrite \eqref{eqTmp8nhgs} once again as
$$ 
e^{-\lp \tau}m_{f} m_{f}^{-1} \left( e^{-\lp \tau}-e^{-\lp_0 \tau}\p_0 \right) + \left( e^{-\lp \tau}-e^{-\lp_0 \tau}\p_0 \right)m_{f}^{-1} m_{f} e^{-\lp_{0} \tau}.
$$
The idea is to \emph{borrow} $r^{-1/4}$ from $e^{-\lp \tau}-e^{-\lp_0 \tau}\p_0$ and \emph{lend} it to $e^{-\lp \tau}$ and $e^{-\lp_{0} \tau}.$  It follows from Theorem~\ref{thmHeatKernBounds} and Equation~\ref{eqSmallHeatKer} that each of the operators
$e^{-\lp \tau}m_{f}, $ $m_{f}^{-1} \left( e^{-\lp \tau}-e^{-\lp_0 \tau}\p_0 \right),$  $\left( e^{-\lp \tau}-e^{-\lp_0 \tau}\p_0 \right)m_{f}^{-1}, $ and $m_{f} e^{-\lp_{0} \tau} $ is Hilbert-Schmidt. Hence $e^{-\lp t}-e^{-\lp_0 t}\p_0$ is of trace-class. For more details, see \cite[Page 259]{Muller3}.
\end{proof}

Now we can apply standard Selberg theory to explicitly evaluate the integral trace  
$$\int_\F \left(  \tr_V K_\Gamma(P,P,t,\chi) - \tr_V k(P,P,t) \right)~dv(P). $$ 

Let $\smat(s)$ be the scattering matrix of $\lp.$ That is the matrix formed from the \emph{constant} terms of the Fourier coefficients of the Eisenstein Series. Let $\phi(s)$ be the determinant of the scattering matrix (see \cite{Friedman1,Friedman2} for more details). 
\begin{thm}\label{thmTraceRelHeat}  
\begin{multline*}
\tr \left( e^{-\lp t}-e^{-\lp_0 t}\p_0   \right)  =\int_\F \left(  \tr_V K_\Gamma(P,P,t,\chi) - \tr_V k(P,P,t) \right)~dv(P)=\\  \sum_{m \in \D} e^{-\lambda_m t} 
-\frac{1}{4\pi} 
\int_{\RR}  \exp\left( -(1+x^2)t \right) \frac{\phi'}{\phi}(ix) \, dx \\
+\frac{1}{4}e^{-t} \tr \smat(0) + \frac{e^{-t}}{\sqrt{4\pi t}} k(\Gamma,\chi) \log Y+ \frac{e^{-t}}{4} k(\Gamma,\chi)
\end{multline*}
\end{thm}
\pf
Since $f(P) \df \tr_V K_\Gamma(P,P,t,\chi) - \tr_V k(P,P,t) \in L^1(\F),$ we can rewrite the integral above as 
\begin{multline}
\lim_{A \ra \infty}\int_{\F_A} \left(  \tr_V K_\Gamma(P,P,t,\chi) - \tr_V k(P,P,t) \right)~dv(P) = \\ \lim_{A \ra \infty} \left( \int_{\F_A}\tr_V K_\Gamma(P,P,t,\chi)~dv(P) -  \int_{\F_A} \tr_V k(P,P,t)~dv(P)    \right). 
\end{multline}
A standard application of the Maa\ss---Selberg relations ( see \cite{Friedman1}, \cite[Page 305]{Elstrodt}, or \cite[Pages 67-70]{Venkov}) gives us 
\begin{multline*}
\int_{\F_A}\tr_V K_\Gamma(P,P,t,\chi)~dv(P) = \\
\log(A)k(\Gamma,\chi)\frac{e^{-t}}{\sqrt{4\pi t}} + \sum_{m \in \D} e^{-\lambda_m t} 
-\frac{1}{4\pi} \int_{\RR}  \exp\left( -(1+x^2)t \right) \frac{\phi'}{\phi}(ix) \, dx 
+\frac{1}{4}e^{-t} \tr \smat(0) + \lto_{A \ra \infty}.
\end{multline*}
A straightforward calculations shows that 
\begin{multline*}
\int_{\F_A}\tr_V k(P,P,t)~dv(P) =  k(\Gamma,\chi) \int_{\F^{Y} \setminus \F_A} \frac{r^2 e^{-t}}{|\pg|\sqrt{4\pi t}}\left(1-\exp\left(\frac{-(\log(r)-\log(Y))^2}{t}  \right)  \right)~dv(P) = \\
k(\Gamma,\chi) \int_{Y}^{A} \int_{\pg}\frac{r^2 e^{-t}}{|\pg|\sqrt{4\pi t}}~\frac{dxdydr}{r^{3}}  - k(\Gamma,\chi) \int_{Y}^{A} \int_{\pg}\frac{r^2 e^{-t}}{|\pg|\sqrt{4\pi t}} \exp\left(\frac{-(\log(r)-\log(Y))^2}{t}  \right) ~\frac{dxdydr}{r^{3}} = \\
\log(A)k(\Gamma,\chi)\frac{e^{-t}}{\sqrt{4\pi t}} - \frac{e^{-t}}{\sqrt{4\pi t}} k(\Gamma,\chi) \log Y -\frac{e^{-t}}{4} k(\Gamma,\chi) + \lto_{A \ra \infty}.
\end{multline*}
The term $\frac{e^{-t}}{4} k(\Gamma,\chi) $ is obtained using a simple $u-$substitution, 
$$u = \frac{\log(r)-\log(Y)}{\sqrt{t}}. $$
\epf

\section{The Relative Zeta-Function and Relative Determinant} \label{secRelZeta}
Now that we have an explicit formula for the trace $$\tr \left( e^{-\lp t}-e^{-\lp_0 t}\p_0   \right)$$ (Theorem~\ref{thmTraceRelHeat}), we can proceed and give an explicit evaluation of the relative zeta-function.

Since  $\sigma(\lp_{0})=\sigma_{\text{ess}}(\lp)= [1,\infty)$ when the representation $\chi$ is singular, 
$$q_{\chi} \df \dim \ker \lp - \dim \ker \lp_{0} = \dim \ker \lp.  $$ 
Following M\"{u}ller we define the relative zeta-function $\zeta(s,\lp,\lp_{0})$ by 
\beq
\zeta(s,\lp,\lp_{0}) = \frac{1}{\Gamma(s)} \int_{0}^{\infty} t^{s-1}\left(\tr \left( e^{-\lp t}-e^{-\lp_0 t}\p_0   \right)-q_{\chi} \right)~dt.
\eeq
Note that in order for the integral above to converge, we need to know the asymptotics of 
$$\tr \left( e^{-\lp t}-e^{-\lp_0 t}\p_0   \right)-q_{\chi} 
$$ 
at both $t=0$ and $t=\infty.$ These asymptotics are given in Lemma~\ref{lemGrowthTheta}. 

Since $\vol(\F) < \infty,$ it follows that \cite[Theorem 3.6.4]{Elstrodt} 
$$
q_{\chi}= \left\{
\begin{array}{cc}
 \dim V  & \mbox{if $\chi$ is trivial}    \\  
    0    & \mbox{else}.        
\end{array}
    \right.
    $$ 
%:ThmRelZetaFunction
\begin{thm}
For $\R(s)>2$ 
\begin{multline*}
\zeta(s,\lp,\lp_{0}) = \sum_{m \in \D} \lambda_m^{-s} 
-\frac{1}{4\pi} 
\int_{\RR}  (1+x^2)^{-s} \frac{\phi'}{\phi}(ix) \, dx 
+\frac{1}{4}\left(\tr\smat(0) + k(\Gamma,\chi) \right) \\ + \frac{k(\Gamma,\chi)}{\sqrt{4\pi}}\frac{\Gamma(s-1/2)}{\Gamma(s)}  \log Y.
\end{multline*}
\end{thm}
\pf
The proof follows from the standard properties of the Mellin transform, Theorem~\ref{thmTraceRelHeat}, and Lemma~\ref{lemGrowthTheta}. Note that $q_{\chi}$ cancels out with any of the terms coming from the \emph{zero} eigenvalues of $\lp.$ 
\epf

\subsection{Asymptotics of the Heat Kernel} The main tool that allows us to study the asymptotic behavior (near $t=0$ and $t=\infty$) of the regularized heat kernel is the Selberg trace formula for the case of a cofinite Kleinian group with finite-dimensional unitary representations \cite{Friedman1,Friedman2}. 

\begin{thm}{(Selberg trace formula)} \label{T:Selberg}
Let $\Gamma $ be a cofinite Kleinian group with one cusp at infinity,  $\chi \in \rep,~h $ be  a holomorphic function on  $ \{ s \in \CC \, | \, |\I(s)| < 2+ \delta \}$ for some $\delta > 0,$ satisfying $ h(1+z^2) = O( 1+|z|^2)^{3/2 - \epsilon}) $ as  $|z| \ra \infty,$ and let
$$ g(x) = \frac{1}{2\pi} \int_{\RR} h(1+t^2)e^{-itx}\,dt. $$  Then
\begin{multline}
\sum _{m \in \D }h(\lambda _{m}) -\frac{1}{4\pi }\int _{\RR }h(1+t^{2})\frac{\phi'}{\phi }(it)\, dt
= \frac{\vol \left( \Gamma \setminus \HH \right)}{4\pi ^{2}}\dim_\CC V \int _{\RR }h(1+t^{2})t^{2}\, dt  \\ + 
 \sum_{ \{R \} \text{\emph{nce}}}\frac{\tr _{V}\chi (R) g(0)\log N(T_{0})}{4|\ren |\sin ^{2}(\frac{\pi k}{m(R)})}  + 
\sum_{\{ T \} \text{\emph{lox}} } \frac{\tr _{V}\chi (T) g(\log N(T))}{\oen |a(T)-a(T)^{-1}|^{2}}\log N(T_{0}) \\ - 
\frac{\tr (\smat (0))h(1)}{4} \\ 
+ \sum_{i=1}^d \frac{\tr \chi(g_i)}{|\mC(g_i)|}\left(  \frac{2 g(0) \log|c_i| }{|1-\epsilon_i^2|^2}+\frac{1}{|1-\epsilon_i^2|^2} \int_0^\infty g(x) \frac{\sinh x}{\cosh x - 1 +\frac{|1-\epsilon_i^2|^2}{2} }~dx \right) \\ 
+ \frac{l_\infty }
{| \Gamma_{\infty}: \Gamma_{\infty}^{\prime} |} 
\left(  g(0)\frac{h(1)}{4} + 
g(0) \left( \frac{ \eta_{\infty}}{2} 
- \gamma \right) - 
\frac{1}{2\pi} 
\int_{\RR} h(1+t^2) \frac{\Gamma'}{\Gamma}(1+it) \,dt   \right)  \\ + 
 \frac{g(0)}{| \Gamma_{\infty}: \Gamma_{\infty}^{\prime} |} \sum_{ l = l_\infty + 1 }^{n} L(\Lambda_{\infty}, \psi_{ l}  ).
\end{multline}
\end{thm}
Here $ \{ \lambda _{m} \}_ {m \in \D } $ are the eigenvalues of $\lp$ counted with multiplicity.  The summation with respect to $\{R \}_\text{nce} $ extends over the finitely many $\Gamma-$conjugacy classes of the non cuspidal elliptic elements (elliptic elements that do not fix a cusp) $R \in \Gamma,$ and for such a class $N(T_0)$ is the minimal norm of a hyperbolic or loxodromic element of the centralizer $\mathcal{C}(R).$  The element $R$ is understood to be a $k-$th power of a primitive non cuspidal elliptic element $R_0 \in \mathcal{C}(R)$ describing a hyperbolic rotation around the fixed axis of $R$ with minimal rotation angle $\frac{2 \pi}{m(R)}.$  Further, $\ren$ is the maximal finite subgroup contained in $  \mathcal{C}(R).$  The summation with respect to $\{ T \}_\text{lox}$ extends over the  $\Gamma-$conjugacy classes of hyperbolic or loxodromic elements of $\Gamma,$ $T_0$ denotes a primitive hyperbolic or loxodromic element associated with $T.$  The element $T$ is conjugate in $\pc$ to the transformation  described by the diagonal matrix with diagonal entries $a(T), a(T)^{-1}$ with $|a(T)| > 1, $ and $N(T) =   |a(T)|^2. $ For $s \in \CC,$ $\smat(s)$ is  a $k(\Gamma,\chi) \times  k(\Gamma,\chi)$ matrix-valued meromorphic function, called the \emph{scattering matrix}  of $\lp,$ and $\phi(s) = \det \smat(s).$ The elements $g_{i}$ are complete representatives for the conjugacy classes of $\{ \gamma \in \gi~|~\gamma~\text{is not parabolic nor the identity element}~\}, $  
$|\mC (g_{i})| $ is the order of the centralizer in $\Gamma $ of the element $g_{i}.$ The numbers $c_{i} \in \CC$ are constants depending on the $g_{i}$ respectively (see \S\ref{sssCE}).
 Finally $\gamma$ is Euler's constant, and $\eta_{\infty}$ is the analogue of Euler's constant for the lattice $\Lambda_{\infty} \subset \RR^{2}.$ The term $L(\Lambda_{\infty}, \psi_{ l}  )$ is defined in \S\ref{sssSiegel}. See \cite{Friedman1,Friedman2,Elstrodt} for more details.

Next, using the Selberg trace formula, we study the regularized heat kernel. We have  
\bl \label{lemGrowthTheta} Let $ \theta(t) \df \tr \left( e^{-\lp t}-e^{-\lp_0 t}\p_0   \right). $ Then there exists constants $a,b,c,d$ so that 
$$\theta(t) =   a t^{-\frac{3}{2}} +b (\log t) t^{-\frac{1}{2}} +c t^{-\frac{1}{2}} + d +O(\sqrt{t}\log t)~\text{as}~t\ra \0^{+},$$
and there exists a positive constant $c>0$ so that $ \theta(t)-q_{\chi} = O(e^{-ct})~\text{as}~t \ra \infty. $
\el
\pf
Our argument is analogous to \cite[Proposition 1]{Efrat1,Efrat2}. 

The asymptotics as $t \ra \infty$ follows immediately from the spectral trace formula (Theorem~\ref{thmTraceRelHeat}).

An application of the Selberg trace formula for the pair of functions 
$$h(z)=\exp(-zt), \quad
g(r)=\frac{\exp(-t)}{\sqrt{4\pi t}}\exp\left(-\frac{r^2}{4t}\right),
$$ 
yields (on applying the left side of the Selberg trace formula) the nontrivial terms from the trace of the heat kernel in Theorem~\ref{thmTraceRelHeat}, namely $$\sum _{m \in \D }e^{-t\lambda _{m}}-\frac{1}{4\pi }\int _{\RR }e^{-t(1+x^{2})}\frac{\phi'}{\phi }(ix)\, dx.  $$
It remains for us to estimate each term on the right as $t \ra 0^{+}.$  

We start with the loxodromic sum 
$$\frac{\exp(-t)}{\sqrt{4\pi t}} \sum_{\{ T \} \text{\emph{lox}} } \frac{\tr _{V}\chi (T) \exp\left(-\frac{(\log N(T))^{2}}{4t}\right)}{\oen |a(T)-a(T)^{-1}|^{2}}\log N(T_{0}). $$ Since $N(T)>1$ for all loxodromic $T,$ the sum decays to zero as $t\ra 0^{+}.$

Next note that $g(0)=\frac{\exp(-t)}{\sqrt{4\pi t}}$ and $h(1)=e^{-t}.$ The finite sum involving the non-cuspidal elliptic terms is easily estimated, 
$$\frac{\exp(-t)}{\sqrt{4\pi t}} \sum_{ \{R \} \text{\emph{nce}}}\frac{\tr _{V}\chi (R) \log N(T_{0})}{4|\ren |\sin ^{2}(\frac{\pi k}{m(R)})} =O(t^{-\frac{1}{2}}), $$ 
\noindent and so are all the other terms with only $g(0)$ or $h(1).$ 

Next we must estimate 
$$
\frac{\exp(-t)}{\sqrt{4\pi t}} \int_0^\infty \exp\left(-\frac{x^2}{4t}\right)   \frac{\sinh x}{\cosh x - 1 +\frac{|1-\epsilon_i^2|^2}{2} }~dx. $$
\noindent An elementary $u-$substitution with $u=\frac{x}{2\sqrt{2}}$ shows that the integral above decays to zero (exponentially fast) as $t \ra 0^{+}.$ 

The next integral can be calculated explicitly: 
$$\frac{\vol \left( \Gamma \setminus \HH \right)}{4\pi ^{2}}\dim_\CC V\int _{\RR }\exp(-t(1+x^{2}))x^{2}\, dx =\frac{\vol \left( \Gamma \setminus \HH \right)}{4\pi ^{2}}(\dim_\CC V)\frac{\sqrt{\pi}e^{-t}}{2 t^{\frac{3}{2}}}.$$
The last integral 
$$\int_{\RR} e^{-(1+x^2)t} \frac{\Gamma'}{\Gamma}(1+ix) \,dx $$ requires some work. To study the integral near $t=0$ we follow \cite[Proposition 1]{Efrat1,Efrat2}. After performing an integration by parts, the integral becomes
$$
-2tie^{-t}\int_{-\infty}^{\infty}xe^{-tx^{2}}\log\Gamma(1+ix)~dx. 
$$ 
Using Sterling's formula 
$$\log\Gamma(1+ix)=(\frac{1}{2}+ix)\log(1+ix)-(1+ix)+\log\sqrt{2\pi}+O(\frac{1}{x})~\text{as}~ x \ra \infty, $$ it follows that for constants $\beta',\gamma',\delta'$
$$\int_{\RR} e^{-(1+x^2)t} \frac{\Gamma'}{\Gamma}(1+ix) \,dx = \beta' \frac{\log t}{\sqrt{t}}+ \frac{\gamma'}{\sqrt{t}}+\delta'+O(\sqrt{t}\log t)~\text{as}~t\ra 0^{+}.$$ \noindent See \cite[Proposition 1]{Efrat1,Efrat2} for more details.
\epf

\subsection{The Regularized Determinant}
In this section we define $\det(\lp - (1-s^{2})).$ 
For $\R(s)>2$ define  $$H(w,s)\df H(w,s,\lp,\lp_{0}) \df \frac{1}{\Gamma(w)} \int_{0}^{\infty} t^{w-1}\left(\tr \left( e^{-\lp t}-e^{-\lp_0 t}\p_0\right)\right)e^{t(1-s^{2})}~dt. $$ 
By applying the Mellin Transform to Theorem~\ref{thmTraceRelHeat} we obtain, for $\R(s)>2,$ 
\begin{multline*}
H(w,s) = \sum_{m \in \D} (\lambda_m-(1-s^{2}))^{-w} 
-\frac{1}{4\pi} 
\int_{\RR}  (x^2+s^{2})^{-w} \frac{\phi'}{\phi}(ix) \, dx 
+s^{-2w}\frac{1}{4}\left(\tr\smat(0) + k(\Gamma,\chi) \right) \\ + s^{-(2w-1)}\frac{k(\Gamma,\chi)}{\sqrt{4\pi}}\frac{\Gamma(w-1/2)}{\Gamma(w)}  \log Y.
\end{multline*}

In order to define $\det(\lp - (1-s^{2})),$ we will need to know that $H(w,s)$ is regular at $w=0.$ 
%:Lemma regularity of H(w,s)
\bl \label{lemRegGrowth} The following hold:
\begin{enumerate}
\item For fixed $s > 2,~H(w,s)$ is regular at $w=0.$
\item \begin{multline*} \frac{\partial H}{\partial w}(0,s)\sim  a (s^{2}-1)^{\frac{3}{2}}+2\sqrt{\pi}b (s^{2}-1)^{\frac{1}{2}}\left( \log(s^{2}-1) +(\gamma +\log(4) - 2) \right) \\ -2\sqrt{\pi}c (s^{2}-1)^{\frac{1}{2}} - d \log(s^{2}-1)~\text{as}~s \ra \infty,  \end{multline*}
\end{enumerate}
where $\gamma$ is Euler's constant, and the constants $a,b,c,d$ are from Lemma~\ref{lemGrowthTheta}.
\el 
\pf
The proof is a standard exercise using the Mellin transform, Lemma~\ref{lemGrowthTheta}, and the following formulas \cite{Grad}:
\begin{align*} 
 &\frac{1}{\Gamma(w)}\int_{0}^{\infty} t^{-\epsilon}e^{t(1-s^{2})}t^{w-1}~dt = \frac{1}{\Gamma(w)} (s^{2}-1)^{w-\epsilon} \Gamma(w-\epsilon),\quad \epsilon=0,\frac{1}{2},\frac{3}{2} \\
 &\frac{1}{\Gamma(w)}\int_{0}^{\infty} \frac{\log t}{\sqrt{t}} e^{t(1-s^{2})}t^{w-1}~dt =\frac{\Gamma(w-\frac{1}{2})}{\Gamma(w)}(s^{2}-1)^{\frac{1}{2}-w}\left( \Psi(w -\frac{1}{2})-\log(s^{2}-1)  \right). \end{align*} 

Here $\Psi(z)$ is the logarithmic derivative of $\Gamma(z).$ Regularity follows from the fact that $\frac{1}{\Gamma(w)}$ vanishes at $w=0.$ See \cite[Prop 2 and 3]{Efrat1,Efrat2} for more details.
\epf

Next, for $\R(s)>2,$ define the \emph{regularized determinant} by
\beq
\det(\lp - (1-s^{2}))=e^{-\frac{\partial H}{\partial w}(0,s)}.
\eeq
Our main result, Theorem~\ref{thmRegDet}, will give the meromorphic continuation to $\R(s) \leq 2.$

\section{Selberg's zeta-function and the regularized determinant} \label{secSelZetaDet} 
\subsection{The Definition of the Selberg zeta-function} \label{ssDefSZF}In this section we define the Selberg zeta-function for cofinite Kleinian groups with finite-dimensional unitary representations. For more details see \cite{Friedman2}.
Suppose $T\in\Gamma$ is loxodromic (we consider hyperbolic elements as loxodromic elements). Then $T$ is conjugate in $\pc$ to a unique element of the form 
$$
D(T)=
\left(\begin{array}{cc}
a(T) & 0\\
0 & a(T)^{-1}
\end{array}\right) $$
such that $a(T)\in\CC$ has $|a(T)|>1$.  Let $N(T)$  denote  the \emph{norm} of $T,$   defined by  $$N(T) \df |a(T)|^{2},$$ and  let   by  $\mC(T) $ denote  the centralizer of $T$ in $\Gamma.$  There exists a (primitive)  loxodromic element $T_0,$ and a finite cyclic elliptic subgroup  $\en$ of order $m(T), $ generated by an element $E_T $   such that 
$$\mC(T) = \langle T_{0} \rangle \times \en. $$
Here $\langle  T_{0} \rangle = \{\, T_{0}^{n} ~ | ~ n \in\ZZ ~ \}. $ Next,
Let $\mathfrak{t}_1,\dots, \mathfrak{t}_n, $ and $\mathtt{t'_1},\dots,  \mathtt{t'_n}$ denote the eigenvalues of $\chi(T_0)$ and $\chi(E_T)$ respectively.   The elliptic element $ E_T$ is conjugate in $\pc$ to an element of the form 
$$\left(\begin{array}{cc}
\zeta(T_0) & 0 \\
0 & \zeta(T_0)^{-1}
\end{array}\right), $$ 
where here $\zeta(T_0)$ is a primitive $2m(T)$-th root of unity.

For $\R(s)>1 $ the Selberg zeta-function $Z(s,\Gamma,\chi)$ is defined by
$$
Z(s,\Gamma,\chi) \df \prod_{ \{T_0 \} \in \mathcal{R}} ~ \prod_{j=1}^{ \dim V} \prod_{  \substack{ l,k \geq 0 \\  c(T,j,l,k)=1   } } \left( 1-\mathfrak{t}_{j} a(T_0)^{-2k} \overline{ a(T_0) ^{-2l}} N(T_0)^{-s - 1}    \right).  
$$
Here the product with respect to $T_0$ extends over a maximal reduced system $\mathcal{R} $ of $\Gamma$-conjugacy classes of primitive loxodromic elements of $\Gamma.$ The system  $\mathcal{R} $ is called reduced if no two of its elements have representatives with the same centralizer\footnote{See \cite{Elstrodt} section 5.4 for more details}.  The function  $c(T,j,l,k)$ is defined by 
$$c(T,j,l,k)= \mathtt{t'_j} \zeta(T_0)^{2l}  \zeta(T_0)^{-2k}.$$
\subsection{The Relationship Between the Selberg Zeta-Function and the Regularized Determinant}
One way to study the Selberg zeta-function is to apply the Selberg trace formula to the pair of functions,  $$h(w)=\frac{1}{s^2+w-1} - \frac{1}{B^2+w-1} ~~\text{and}~ $$ 
$$ g(x) = \frac{1}{2s}e^{-s|x|} - \frac{1}{2B}e^{-B|x|}, $$ where $1 < \R(s) < \R(B).$  
Let $Z(s)\df Z(s,\Gamma,\chi)$ be the Selberg zeta-function under Assumption~\ref{asOne}. We have \cite{Friedman2}:
\begin{lem}
\begin{multline} \label{eqLogDer}
\frac{1}{2s} \frac{Z^\prime}{Z}(s) - \frac{1}{2B} \frac{Z^\prime}{Z}(B)  
=\frac{1}{2s} \sum_{ \{ T \}\LOX}  \frac{   \tr (\chi(T)) \log N(T_{0})}{m(T)|a(T)-a(T)^{-1}|^{2}}N(T)^{-s} \\
 -\frac{1}{2B}\sum_{ \{ T \}\LOX} \frac{ \tr (\chi(T))  \log N(T_{0})}{m(T)|a(T)-a(T)^{-1}|^{2}}N(T)^{-B} \\
= \sum_{n \in D} \left(\frac{1}{\lambda_{n}-(1-s^{2})} - \frac{1}{\lambda_{n}-(1-B^{2})}  \right) 
- \frac{1}{4 \pi} \int_\RR  \left(\frac{1}{s^2 +x^2} - \frac{1}{B^2 + x^2}  \right) \frac{\phi^\prime}{\phi}(i x)~dx  \\
+ \frac{l_\infty}{2 \pi [\gi:\gip]} \int_\RR  \left(\frac{1}{s^2 +x^2} - \frac{1}{B^2 + x^2}  \right) \frac{\Gamma^\prime}{\Gamma}(1 + ix)~dx 
+ \frac{\tr \smat(0)}{4s^2} - \frac{\tr \smat(0)}{4B^2} \\
- \frac{l_\infty}{4[\gi:\gip]s^2} +\frac{l_\infty}{4[\gi:\gip]B^2}  \\
- \sum_{i=1}^l \frac{\tr \chi(g_i) }{|C(g_i)||1-\epsilon_i^2|^2} \int_0^\infty  \left( \frac{e^{-sx}}{2s} - \frac{e^{-Bx}}{2B} \right)   \frac{\sinh x}{\cosh x -1 +\frac{|1-\epsilon_i^2|^2}{2} }~dx 
   \\ - \left( \frac{1}{2s}-\frac{1}{2B}\right) \sum_{ \{R \} \text{\emph{nce}}}\frac{\tr _{V}\chi (R) \log N(T_{0})}{4|\ren |\sin ^{2}(\frac{\pi k}{m(R)})}
+ \frac{\vol \left( \Gamma \setminus \HH \right)\dim V  }{4\pi}(s-B)
 \\ - \left( \frac{1}{2s}-\frac{1}{2B}\right) \sum_{i=1}^l \frac{2 \tr \chi(g_i) \log|c_i| }{|C(g_i)||1-\epsilon_i^2|^2} 
 - \left( \frac{1}{2s}-\frac{1}{2B}\right) \frac{1}{[\gi:\gip]} \left(  l_\infty \left( \frac{ \eta_{\infty}}{2} 
- \gamma \right) +  \sum_{ l = l_\infty + 1 }^{n} L(\Lambda_{\infty}, \psi_{ l}  )   \right).
\end{multline}
\end{lem}

Recall that \begin{multline*}
H(w,s) = \sum_{m \in \D} ((\lambda_m-1)+s^{2})^{-w} 
-\frac{1}{4\pi} 
\int_{\RR}  (x^2+s^{2})^{-w} \frac{\phi'}{\phi}(ix) \, dx 
+s^{-2w}\frac{1}{4}\left(\tr\smat(0) + k(\Gamma,\chi) \right) \\ + s^{(1-2w)}\frac{k(\Gamma,\chi)}{\sqrt{4\pi}}\frac{\Gamma(w-1/2)}{\Gamma(w)}  \log Y.
\end{multline*}

Applying the following elementary equations: 
$$-\frac{d}{ds}\left(\frac{1}{2s}\frac{d}{ds}\left(\frac{\partial}{\partial w} \left. (u+s^{2})^{-w}\right|_{w=0}  \right)  \right) = \left. - \frac{\partial}{\partial w}\frac{d}{ds}\left(\frac{1}{2s}\frac{d}{ds}  (u+s^{2})^{-w}    \right) \right|_{w=0}  = \frac{-2s}{(u+s^{2})^{2}}, $$
%:H(w,s) differentiated
$$-\frac{d}{ds}\left(\frac{1}{2s}\frac{d}{ds}\left(\frac{\partial}{\partial w} \left.  s^{(1-2w)}\frac{\Gamma(w-1/2)}{\Gamma(w)} \right|_{w=0} \right)  \right) = -\frac{\sqrt{\pi}}{s^{2}},  $$
to $H(w,s),$ we see that
\begin{multline} \label{eqDiffZeta}
-\frac{d}{ds}\left(\frac{1}{2s}\frac{d}{ds}\left(\frac{\partial}{\partial w} \left. H(w,s) \right|_{w=0}  \right)  \right) = \left. - \frac{\partial}{\partial w}\frac{d}{ds}\left(\frac{1}{2s}\frac{d}{ds}  H(w,s)    \right) \right|_{w=0}  = \sum_{m \in \D} \frac{-2s}{((\lambda_m-1)+s^{2})^{2}} 
\\ -\frac{1}{4\pi} \int_{\RR}  \frac{-2s}{(x^2+s^{2})^{2}} \frac{\phi'}{\phi}(ix) \, dx 
  - \frac{k(\Gamma,\chi)}{2s^{2}}\log Y.
\end{multline}

\textbf{Caution.} Differentiation through the sum and integral is justified by regularity at $w=0$ (Lemma~\ref{lemRegGrowth}). Instead of differentiating first with respect to $w$ at $w=0,$ we switch the order of differentiation, differentiate with respect to $s,$ and restrict $w$ so that $w > 2$ (where the sum and integral converge uniformly). Finally we differentiate with respect to $w,$ and using analytic continuation (and uniqueness of analytic continuation), we obtain \eqref{eqDiffZeta}. 
 
%:logder of zeta
For $\R(s)>0$ set 
\beq \label{eqOmEll} \Omega(s)=-\sum_{i=1}^l \frac{\tr \chi(g_i) }{|C(g_i)||1-\epsilon_i^2|^2}\int_0^\infty e^{-sx}  \frac{\sinh x}{x} \left( \frac{1}{\cosh x -1 +\frac{|1-\epsilon_{i}^2|^2}{2} }\right)~dx, ~\text{for}~\epsilon_{i} \neq 1, \eeq (see \S\ref{sssCE} for the definition of $\epsilon_{i}$ and $g_{i}).$
Using the equation $$\frac{1}{\pi}\int_\RR  \frac{1}{(s^2 +x^2)} \frac{\Gamma^\prime}{\Gamma}(1 + ix)~dx =\frac{1}{s}\frac{\Gamma^\prime}{\Gamma}(1 + s), $$
it follows that  
\begin{multline} \label
{eqDerLogZeta}
\frac{d}{ds}\left(\frac{1}{2s} \frac{Z^\prime}{Z}(s) \right)=  \sum_{n \in D} \frac{-2s}{(\lambda_{n}-1+s^{2})^{2}} 
- \frac{1}{4 \pi} \int_\RR  \frac{-2s}{(s^2 +x^2)^{2}} \frac{\phi^\prime}{\phi}(i x)~dx  \\
+ \frac{d}{ds}\frac{l_\infty}{2s[\gi:\gip]}\frac{\Gamma'}{\Gamma}(s+1) 
-\frac{1}{2s^{3}}\left(\tr \smat(0) 
- \frac{l_\infty}{[\gi:\gip]}\right)   \\
- \frac{d}{ds} \frac{1}{2s}(\frac{d}{ds}\Omega(s))  +  \frac{1}{2s^{2}} \sum_{ \{R \} \text{\emph{nce}}}\frac{\tr _{V}\chi (R) \log N(T_{0})}{4|\ren |\sin ^{2}(\frac{\pi k}{m(R)})}
 \\ +\frac{1}{2s^{2}} \left[ \sum_{i=1}^l \frac{2 \tr \chi(g_i) \log|c_i| }{|C(g_i)||1-\epsilon_i^2|^2} 
 + \frac{1}{[\gi:\gip]} \left(  l_\infty \left( \frac{ \eta_{\infty}}{2} 
- \gamma \right) +  \sum_{ l = l_\infty + 1 }^{n} L(\Lambda_{\infty}, \psi_{ l}  )   \right) \right] 
\\ + \frac{\vol \left( \Gamma \setminus \HH \right)\dim V  }{4\pi}.
\end{multline}

Simplifying, and recalling the definition of $\det(\lp-(1-s^{2}))$ we arrive at 
\begin{multline*}
\frac{d}{ds}\left(\frac{1}{2s} \frac{Z^\prime}{Z}(s) \right)= \frac{d}{ds}\left(\frac{1}{2s}\frac{d}{ds}\left( \log(\det(\lp-(1-s^{2})) \right)  \right) + \frac{k(\Gamma,\chi)}{2s^{2}}\log Y - \frac{d}{ds}\frac{1}{2s}\left(\frac{l_\infty}{[\gi:\gip]}\frac{\Gamma'}{\Gamma}(s+1) \right) \\
-\frac{d}{ds}\left(\frac{1}{2s} \frac{d}{ds}\Omega(s)\right) +\frac{1}{2s^{2}}C_{1}-\frac{1}{2s^{3}}C_{2}+C_{3},
\end{multline*}
where $C_{1},C_{2},C_{3}$ are easily read off from \eqref{eqDerLogZeta}.
Next, integrating twice, we obtain:
\begin{multline*}
\log Z(s) = \log\det \left(\lp-(1-s^{2}) \right) -(s) k(\Gamma,\chi)\log Y - \frac{l_\infty}{[\gi:\gip]}\log \Gamma(s+1)  \\
-\Omega(s) -sC_{1}+\frac{C_{2}}{2}\log s+\frac{2}{3}C_{3}s^{3}+D_{1}+s^{2}D_{2},
\end{multline*}
where $D_{1},D_{2}$ are constants of integration. They can be determined by letting $s \ra \infty,$ and applying Lemma~\ref{lemRegGrowth}. More specifically, Lemma~\ref{lemRegGrowth} tells us the asymptotic growth of $\log\det \left(\lp-(1-s^{2}) \right) $ as $s \ra \infty.$ Noting that $$\lim_{s\ra \infty}\log Z(s) = 0, \quad \lim_{s\ra \infty} \Omega(s)=0,$$ and applying Sterling's formula we obtain:
$$D_{1}=\gamma +\log(4) - 2 + \frac{l_\infty}{[\gi:\gip]}\log\sqrt{2\pi},\quad D_{2}=0.  $$
We have proved 
\begin{thm} \label{thmRegDet} For $\R(s) > 2,$ 
\begin{multline*}
\log\det \left(\lp-(1-s^{2}) \right) = \log Z(s,\Gamma,\chi) + s \left( k(\Gamma,\chi)\log(Y) +C_{1}\right) \\ +\frac{l_\infty}{[\gi:\gip]}\log \Gamma(s+1) 
+ \Omega(s)  - \frac{C_{2}}{2}\log s-\frac{2}{3}C_{3}s^{3}-D_{1},
\end{multline*}
where 
\begin{multline*}
C_{1}= \sum_{ \{R \} \text{\emph{nce}}}\frac{\tr _{V}\chi (R) \log N(T_{0})}{4|\ren |\sin ^{2}(\frac{\pi k}{m(R)})} \\ +\left[ \sum_{i=1}^l \frac{2 \tr \chi(g_i) \log|c_i| }{|C(g_i)||1-\epsilon_i^2|^2} 
 + \frac{1}{[\gi:\gip]} \left(  l_\infty \left( \frac{ \eta_{\infty}}{2} 
- \gamma \right) +  \sum_{ l = l_\infty + 1 }^{n} L(\Lambda_{\infty}, \psi_{ l}  )   \right) \right],   
\end{multline*}
$$C_{2}=\left(\tr \smat(0) 
- \frac{l_\infty}{[\gi:\gip]}\right), $$
$$C_{3}=\frac{\vol \left( \Gamma \setminus \HH \right)\dim V  }{4\pi},    $$ and
$$D_{1}=\gamma +\log(4) - 2 + \frac{l_\infty}{[\gi:\gip]}\log\sqrt{2\pi}  $$ 
\end{thm}

\begin{cor}
Let $\Gamma$ be torsion-free with one cusp at $\infty,$ and let $\chi$ be a regular character (a one-dimensional unitary representation). Then 
$$\det \left(\lp-(1-s^{2}) \right) = Z(s,\Gamma,\chi)\exp\left(-s^{3}\frac{\vol \left( \Gamma \setminus \HH \right)}{6\pi}+sL(\Lambda_{\infty},\psi) \right). $$
\end{cor}

\begin{cor}
Let $\Gamma$ be cocompact, and let $\chi$ be a regular character. Then 
$$\det \left(\lp-(1-s^{2}) \right) = Z(s,\Gamma,\chi)\exp\left(-s^{3}\frac{\vol \left( \Gamma \setminus \HH \right)}{6\pi}+sC \right), $$
where $$C=\sum_{ \{R \} \text{\emph{nce}}}\frac{\tr _{V}\chi (R) \log N(T_{0})}{4|\ren |\sin ^{2}(\frac{\pi k}{m(R)})}.$$
\end{cor}

\appendix

\section{Proof of Theorem~\ref{thmHeatKernBounds}}    \label{appHeatKernBounds}
In this section we prove Theorem~\ref{thmHeatKernBounds}, 
That is for $P=z+rj, P'=z'+r'j \in \F^{Y},$ we show that  $$K_\Gamma(P,P',t,\chi) =\prV rr'\frac{\exp(-t)}{|\pg| \sqrt{4\pi t}}\exp\left(-\frac{\log^{2}(r/r')}{4t}\right)+O(1).$$ 
As usual we are under Assumption~\ref{asOne}.

The first step is to split up  
$K_\Gamma(P,P',t,\chi)$ as 
\begin{multline*} K_\Gamma(P,P',t,\chi) =\sum_{\gamma \in \Gamma} \chi(\gamma) k_t(\delta(P,\gamma P'))=  \sum_{\gamma \in \gi} \chi(\gamma) k_t(\delta(P,\gamma P'))  + \sum_{\gamma \in \Gamma \setminus \gi} \chi(\gamma) k_t(\delta(P,\gamma P')). \end{multline*} 
From \cite[Equation 4.5.9 or Lemma 6.4.2]{Elstrodt}, it follows that 
$$  \sum_{\gamma \in \Gamma \setminus \gi} |k_t(\delta(P,\gamma P'))| = O(1)~\text{as}~r\ra \infty. $$
Hence 
$$\sum_{\gamma \in \Gamma \setminus \gi} \chi(\gamma) k_t(\delta(P,\gamma P')) = O(1) \quad \text{as $r \ra \infty.$}$$
It remains to estimate $$f_\infty(P,P') \df \sum_{\gamma \in \gi} \chi(\gamma) k_t(\delta(P,\gamma P')).  $$ 

The subgroup $\Gamma_\infty$ is not an abelian group. So in general we can not diagonalize the unitary representation, $\chi$ restricted to $\Gamma_\infty,$ into unitary characters. However, we have the following lemma which is almost as good as diagonalizing $\chi$ (see \cite[Lemma 2.4]{Friedman2}).
\begin{lem}
\label{L:mylemma} 
Let $\infty$ be the one cusp of $\Gamma$ (Assumption~\ref{asOne}). Then there exist $E,R,S \in \gi $ with the following properties:
\begin{enumerate}
\item \label{partOneMylemma} $\gi=\{\, E^{k}R^{i}S^{j}\,|\,0\leq k<m,\, i,j\in\ZZ\,\}.$
Here $R,S$ are parabolic elements with $ R(P) = P+1$ and $S(P) = P+\tau$ (here $\Lambda_\infty = \ZZ \oplus \ZZ \tau$) for all $P \in \HH,$  and  $E$ is elliptic of order $m.$ 

\item  $\gip =\{\, R^{i}S^{j}\,|\, i,j\in\ZZ\,\}.$

\item The elements $R$ and $S$ commute, but the group $\gi$ is not abelian when $m>1$. 

\item If in addition, $m > 1, $  
then  $\chi(E)$ maps $V_{\infty}^\prime $ onto itself. Furthermore, there exists
a basis of $V_{\infty}^\prime$ so that $\chi(E)|_{V_{\infty}^\prime}$ is diagonal. 
\end{enumerate}
\end{lem}
\noindent The notation used in the lemma above is explained in \S\ref{secNotation}.

Next we split $f_\infty$ into two sums,
$$
f_\infty(P,P') =  \sum_{\gamma \in \gi} \prV \chi(\gamma) k_t(\delta(P,\gamma P')) + \sum_{\gamma \in \gi} \left( I_V - \prV\right) \chi(\gamma) k_t(\delta(P,\gamma P')).
$$

\bl[Poisson Summation Formula]
Let $f:\RR^2 \mapsto \CC$ be a Schwartz-class function, and let $\Lambda$ be a two-dimensional lattice in $\RR^2=\CC.$ Then 
$$\sum_{\omega \in \Lambda} f(\omega) = \frac{1}{|\Lambda|} \sum_{\omega \in \Lambda^0} \widehat{f}(\omega).  $$ 
Here $\Lambda^0$ is the dual lattice to $\Lambda,$ and $\widehat{f}$ is the Fourier transform of $f,$
$$
\widehat{f}(z) = \int_{\RR^2} f(u)e^{-2\pi i \left< u, z \right>}~du,  
$$
where  $\left< \cdot, \cdot \right>$ is the standard real inner product on $\CC = \RR^2.$
\el

\bl \label{lemSingKerBound}
$$ \sum_{\gamma \in \gi} \prV \chi(\gamma) k_t(\delta(P,\gamma P')) = \prV   rr'\frac{\exp(-t)}{|\pg| \sqrt{4\pi t}}\exp\left(-\frac{\log^{2}(r/r')}{4t}\right) + O(1).  $$
\el
\pf
It follows from the definition of singularity that for $\gamma \in \Gamma_\infty,$ $\prV \chi(\gamma) = \prV.$ Thus 
$$\sum_{\gamma \in \gi} \prV \chi(\gamma) k_t(\delta(P,\gamma P')) = \prV \sum_{\gamma \in \gi}  k_t(\delta(P,\gamma P')).    $$ 
Using  Part~\ref{partOneMylemma} of Lemma~\ref{L:mylemma}, we can rewrite the above sum as
$$\prV \sum_{k=0}^{m}\sum_{i,j \in \ZZ^2}k_t(\delta(P,E^k R^i S^j P'))=\prV \sum_{k=0}^{m}\sum_{i,j \in \ZZ^2}k_t(\delta(E^{-k}P, R^i S^j P')).  
$$
The last equality follows because $\delta$ is a point-pair invariant. Once again applying Lemma~\ref{L:mylemma}, and the definition of $\delta$ we can write the above sum as 
$$ \prV \sum_{k=0}^{m} \sum_{\omega \in \Lambda_\infty} f_k(\omega),
$$
where $$f_k(u)=k_t\left( \frac{|z-z(E^k P')+u|^2+r^2 + r'^{2}  }{2rr'} \right). $$ We have used\footnote{We abuse notation here and let $r,z$ represent both coordinates in $\HH$ and coordinate functions.} the fact that $r(E^k P) = r(P) = r.$ Next we apply the Poisson summation to obtain 
$$ \sum_{\omega \in \Lambda_\infty} f_k(\omega) = |\Lambda_\infty|^{-1} \hat{f_k}(0) + |\Lambda_\infty|^{-1} \sum_{\substack{ \omega \in \Lambda_\infty \\ \omega \neq 0 }} \hat{f_k}(\omega). $$ A straightforward computation \cite[Lemma 3.5.5]{Elstrodt}
shows\footnote{\cite[Lemma 6.4.2]{Elstrodt} seems to have an extra factor of ``2.''  }  that 
$$
\hat{f_k}(0) = rr'g\left( \log\left( \frac{r}{r'} \right)  \right)$$
Where 
$$
g(x) = \frac{\exp(-t)}{\sqrt{4\pi t}}\exp\left(-\frac{x^2}{4t}\right).
$$
Noting that $\hat{f_k}(0)$ is independent of $k,$ and that $$\frac{m}{|\Lambda_\infty|} = \frac{|E|}{|\Lambda_\infty|} = \frac{[\gi:\gip]}{|\Lambda_\infty|} = \frac{1}{|\pg|}$$ we recover the leading term of the lemma. To conclude we show that $$\sum_{\substack{ \omega \in \Lambda_\infty \\ \omega \neq 0 }} \hat{f_k}(\omega)=O(1). $$ Since $f_k$ is smooth, $\hat{f_k}$ is of rapid decay, hence \cite[Lemma 6.4.2]{Elstrodt}
\beq \label{eqDecayUpperTerms}
\hat{f_k}(v)=O( (rr')^{-N} |v|^{-2N})~\text{for any $N>0$}.
\eeq
The lemma now follows.
\epf
Equation~\ref{eqDecayUpperTerms} will be used to show that  $$\sum_{\gamma \in \gi} \left( I_V - \prV\right) \chi(\gamma) k_t(\delta(P,\gamma P')) = O(1) ~\text{as}~r \ra \infty.$$  In order to proceed, we need to understand $I_V - \prV$ as projection operator on the subspace $V_\infty^\perp.$ First, decompose $V$ into 
$$V = V_\infty \oplus V_\infty^\prime \oplus V_\infty^b. $$ 
By Lemma~\ref{L:mylemma} the unitary representation $\chi$ restricted to $\Gamma_\infty$ can be diagonalized in a block matrix form with respect to the decomposition $V = V_\infty \oplus V_\infty^\prime \oplus V_\infty^b. $ Hence we can write $I_V - \prV = \mathbf{P}_a + \mathbf{P}_b,$ where $\mathbf{P}_a$ is the orthogonal projection onto $V_\infty^\prime$ and $\mathbf{P}_b$ is the orthogonal projection onto $V_\infty^b.$   
\bl
$$\sum_{\gamma \in \gi} \left( I_V - \prV\right) \chi(\gamma) k_t(\delta(P,\gamma P')) = O(1).$$
\el
\pf
We first show that $$\sum_{\gamma \in \gi} \mathbf{P}_a \chi(\gamma) k_t(\delta(P,\gamma P')) = O(1).$$ The proof is almost identical to the proof of Lemma~\ref{lemSingKerBound} except that the term corresponding to $\omega=0,$ is zero. Indeed, 
$$\mathbf{P}_a \sum_{k=0}^{m}\sum_{i,j \in \ZZ^2}\chi(E^k R^i S^j) k_t(\delta(P,E^k R^i S^j P')) = \left( \sum_{k=0}^{m} \chi(E)^k \right) \mathbf{P}_a\sum_{i,j \in \ZZ^2}k_t(\delta(P,E^k R^i S^j P')),$$ where $\chi(E)$ is a diagonal matrix with each element on the diagonal a finite root of unity not equal to one. Since the order of each root of unity divides $m$ we must have $$\sum_{k=0}^{m} \chi(E)^k = 0. $$ Hence the ``constant term'' (the term corresponding to $\omega=0$ cancels out).

Next, to show that $$\sum_{\gamma \in \gi} \mathbf{P}_b \chi(\gamma) k_t(\delta(P,\gamma P'))$$ is bounded it suffices to estimate the lattice sum
\beq \label{eqTmp8431d} \mathbf{P}_b\sum_{l,j \in \ZZ^2}k_t(\delta(P,R^l S^j P')).  \eeq
Since $R$ and $S$ commute, we can diagonalize $\chi(\gip)$ restricted to $V_\infty^b.$ Hence we can assume that $\chi$ is a lattice character of the form 
$$\chi(R^lS^n) = \exp(2\pi i (l \theta_R+ n \theta_S)). $$ Now we can rewrite \eqref{eqTmp8431d} as 
$$ 
\sum_{l,n\in\ZZ^{2}}\exp(2\pi i (l \theta_R+ n \theta_S))k_t(\delta(z+rj,l+n\tau (z'+r'j) ))
$$
By unraveling the definition of $\mathbf{P}_{b},$ it follows that at least one of $\theta_R,\theta_S$ is not an integer. By applying the Poison summation formula  to the function 
$$ f_{1}(w,v) = \exp(2\pi i (w \theta_R+ v \theta_S))k_t(\delta(z+rj,w+v\tau (z'+r'j) )), $$ we obtain
$$\sum_{l,n\in \ZZ^{2}}f_{1}(l,n) = \sum_{l,n\in \ZZ^{2}}\widehat{f_{1}}(l,n).$$ However the exponential factor $\exp(2\pi i (l \theta_R+ n \theta_S))$ shifts the Fourier transform of the function $f_{2}(w,v)=k_t(\delta(z+rj,w+v\tau (z'+r'j) )),$ and wipes out the unbounded ``constant term'' $\widehat{f_{2}}(0,0).$ In other words, if we applied the Poisson summation formula to $f_{2},$ we would see, as we did with $f_{k},$ that 
$$\sum_{l,n\in \ZZ^{2}}f_{2}(l,n) = \widehat{f_{2}}(0,0) + \sum_{\substack{l,n\in \ZZ^{2} \\ (l,n)\neq(0,0)}}\widehat{f_{2}}(l,n).$$ The latter sum decays, while the first \emph{term} would not. The effect of multiplying $\exp(2\pi i (w \theta_R+ v \theta_S))$ is to shift the sum away from the integers. That is
$$\sum_{l,n\in \ZZ^{2}}\widehat{f_{1}}(l,n) = \sum_{l,n\in \ZZ^{2}}f_{2}(l+\alpha,n+\beta),$$ where
$$(0,0) \notin \ZZ^{2}+(\alpha,\beta).$$ Hence we can apply \eqref{eqDecayUpperTerms} to conclude the lemma.  
\epf

\bibliography{reg}
\bibliographystyle{amsalpha}
\end{document}